\documentclass[11pt]{amsart}

\usepackage{AFBoix2015}

\begin{document}

\author[A.\,F.\,Boix]{Alberto F.\,Boix}
\address{Department of Mathematics, Universitat Polit\`ecnica de Catalunya BarcelonaTech, Av. Eduard Maristany 16, 08019, Barcelona, Spain.}
\email{alberto.fernandez.boix@upc.edu}

\author[S. Zarzuela]{Santiago Zarzuela}
\address{Departament de Matem\`atiques i Inform\`atica, Universitat de Barcelona, Gran Via de les Corts Catalanes 585, 08007, Barcelona, SPAIN}
\email{szarzuela@ub.edu}

\title[On deficiency modules]{Regularity of deficiency modules through spectral sequences}



\begin{abstract}
The main goal of this paper is to obtain upper bounds for the regularity of graded deficiency modules in the spirit of the one obtained by Kummini--Murai in the monomial case building upon the spectral sequence formalism developed by \`Alvarez Montaner, Boix and Zarzuela. This spectral sequence formalism allows us not only to recover Kummini--Murai's upper bound for monomial ideals, but also to extend it for other types of rings, which include toric face rings and some binomial edge rings, producing to the best of our knowledge new upper bounds for the regularity of graded deficiency modules of this type of rings.
\end{abstract}

\maketitle

\section*{Introduction}
Let $\mathbb{K}$ be a field, let $R=\mathbb{K}[x_1,\ldots, x_n]$ be a polynomial ring in $n$ variables over $\mathbb{K},$ that we assume it has the standard grading (this just means that elements of $\mathbb{K}$ have degree zero, and variables have degree one). In this setting, given a finitely generated, graded $R$--module $M$ we consider Schenzel's graded deficiency modules
\[
K^j (M):=\hbox{}^* \Ext_R^{n-j} (M,R(-n))
\]
as a way of measuring how far $M$ is of being a Cohen--Macaulay module. In this setting, it is of interest to understand how big can be the Castelnuovo--Mumford regularity of these modules; more precisely, we want to calculate upper bounds for this regularity.

This is a topic that has been considered in several research papers, where general upper bounds has been obtained (see, for instance, \cite{HoaHyry2006,BrodmannLinhSeiler2013,ChardinHaHoa2011}). However, most of these bounds are far from being explicit.

Now, we assume that $M=R/I$ where $I$ is a monomial ideal; in this setting, Kummini and Murai \cite{KumminiMurai2011} obtained the following nice and explicit upper bound:
\[
\reg (K^j (R/I))\leq j.
\]
The proof of this result deeply relies on the so--called \textit{Stanley filtrations} used for instance by Maclagan and Smith \cite{MaclaganSmith2005} to obtain upper bounds for the so--called multigraded regularity for toric varieties that they introduce in \cite{MaclaganSmith2004}; however, Stanley filtrations require working with multigradings, a technique which is almost exclusive for monomial ideals.

In other different direction, Ha and Hoa \cite[Proposition 2.3]{HaHoa2008} showed that if $I$ is an arbitrary graded ideal such that $R/I$ is Cohen--Macaulay, then
\[
\reg (K^{\dim (R/I)}(R/I))=\dim (R/I).
\]
For us, one of the interests for studying the regularity of deficiency modules is given by \cite{DaoMaVarbaro}; indeed, assume now that $I$ is a graded ideal inside $R,$ and let $(h_a (R/I),\ldots, h_s (R/I))$ be the $h$--vector of $R/I.$ By a celebrated result by Murai and Terai \cite[Theorem 1.4]{MuraiTeraihvectors}, it is known that, if there exists an integer $r\geq 0$ such that
\[
\reg (K^j (R/I))\leq j-r
\]
for all $j<\dim (R/I),$ then $h_j (R/I)\geq 0$ for all $j\leq r.$ This result was extended by Dao, Ma and Varbaro in \cite[Theorem 2.8]{DaoMaVarbaro} replacing $R/I$ by an arbitrary graded, finitely generated $R$--module $M.$ The interested reader can consult for instance Holmes work \cite[Corollary 5.5]{HolmesgeneralizedSerrecondition} for a setting where one can obtain this type of upper bounds for the regularity.

The main goal of this paper is to produce upper bounds in the spirit of the one obtained by Kummini--Murai in the monomial case building upon the spectral sequence formalism developed in \cite{AlvarezBoixZarzuelaspectralsequences}. This spectral sequence formalism allows us not only to recover Kummini--Murai's upper bound for monomial ideals, but also to extend it for other types of rings.

Now, we provide a more detailed overview of the contents of this paper for the convenience of the reader. After reviewing in Section \ref{section on preliminaries on spectral sequences} some preliminary results obtained in \cite{AlvarezBoixZarzuelaspectralsequences} that we plan to use along this paper, in Section \ref{section on ultrametric functions} we introduce what we call ultrametric functions, that can be regarded as a generalization of the notion of regularity. We do so mainly in order to illustrate that the formalism used throughout this paper can be applied not only for the regularity, but also for other functions that behave, roughly speaking, like regularity on short exact sequences.

Section \ref{section of main results} contains the main results of this paper; here, we want to single out Theorem \ref{the expected bound when P is a C.M. poset}, where we establish our main upper bound for the regularity of graded deficiency modules. Theorem \ref{the expected bound when P is a C.M. poset} recovers and extends \cite{KumminiMurai2011} beyond the monomial case.

In the remainder sections we specialize Theorem \ref{the expected bound when P is a C.M. poset} to some concrete settings; more precisely, in Section \ref{section: the case of toric face rings} we specialize to the case of toric face rings, obtaining, to the best of our knowledge, a new upper bound for the regularity of deficiency modules of a toric face ring (see Theorem \ref{deficiency in case of toric face rings}). On the other hand, in Section \ref{section: the case of binomial edge ideals} we also obtain similar upper bounds for some binomial edge rings (see Theorem \ref{deficiency in case of binomial edge ideals}). Finally, in Section \ref{section of examples} we illustrate our main results by means of examples.


\section{Preliminaries on some spectral sequences}\label{section on preliminaries on spectral sequences}

\subsection{Construction of homological spectral sequences}
The purpose of this section is to collect some results borrowed from \cite{AlvarezBoixZarzuelaspectralsequences} that we plan to use along this paper.
The setup we need for the appropriate  spectral sequence is the following:

\begin{cons}\thlabel{otra construccion mas de lo mismo}
Let $\mathbb{K}$ be a field, let $A=\mathbb{K}[x_1,\ldots ,x_n]$ be the polynomial ring in $n$ variables over $\mathbb{K}$ having the standard grading, and let $\hbox{}^*\modfg_A$ be the category of fintely generated graded $A$--modules. Moreover, let $P$ be any finite poset, let $\Inv (P,\hbox{}^*\modfg_A)$ be the category of inverse systems on $P$ valued on $\hbox{}^*\modfg_A$, and let $\Dir (P,\hbox{}^*\modfg_A)$ be the category of direct systems on $P$ valued on $\hbox{}^*\modfg_A$. We define a bivariate functor (namely, $\mathcal{T}$) in the following manner:
\begin{align*}
& \xymatrix{\Inv (P,\hbox{}^*\modfg_A)\times \hbox{}^*\modfg_A\ar[r]^-{\mathcal{T}}& \Dir (P,\hbox{}^*\modfg_A)}\\ & (G=(G_p)_{p\in P},N)\longmapsto
\mathcal{T}(G,N):=(\hbox{}^*\Hom_A (G_p,N))_{p\in P}.
\end{align*}
Observe that, if we fix the first argument, then $\mathcal{T}(G,-)$ is a covariant, left exact functor, while if we fix the second argument, then $\mathcal{T}(-,N)$ becomes a contravariant, left exact functor. On the other hand, given $M$ a finitely generated graded $A$--module, we denote by $M_{\leq q}$ the following inverse system.
\[
(M_{\leq q})_p=\begin{cases}
    M,\text{ if }p\leq q,\\
    0,\text{ otherwise.}
\end{cases}
\]
Next, given $q\in P$, and given once again $M$ a finitely generated graded $A$--module, we denote by $M_q$ the following direct system: for any $p\in P,$
\[
(M_q)_p=\begin{cases}
M,\text{ if }p=q,\\
0,\text{ otherwise.}
\end{cases}
\]
We also recall that, given $t\in\mathbb{Z}$, $M(t)$ denotes the $t$--shift of $M$ as graded $A$--module. This just means that, if $M=\oplus_{n\in\mathbb{Z}} M_n$, then $(M(t))_n=M_{n+t}$. Finally, given $V\in\Inv (P,\hbox{}^*\modfg_A)$, we say that $V$ is acyclic with respect to the limit functor provided
\[
\left(\mathbb{R}^j \lim_{p\in P}\right) (V_p)=0\text{ for any }j\geq 1.
\]
\end{cons}
Our plan now is to prove the following technical statement, that will be crucial for us along the main result of this section.

\begin{lm}\thlabel{key lemma}
Preserving the notations and assumptions of \thref{otra construccion mas de lo mismo}, the following assertions hold.

\begin{enumerate}[(i)]

\item Let $\mathbb{R}\mathcal{T}$ be the right derived functors of $\mathcal{T}$ regarded as bivariate functor \cite[Chapter V, \S 3]{CartanEilenberg1956}. Moreover, let $\mathbb{R}_{(1)}\mathcal{T}$ and $\mathbb{R}_{(2)}\mathcal{T}$ be the right partial derived functors of $\mathcal{T}$ \cite[Chapter V,\S 8]{CartanEilenberg1956}, where the subscript means that we fix either the first or the second argument of $\mathcal{T}$. Then, we have that $\mathbb{R}_{(1)}\mathcal{T}\cong\mathbb{R}\mathcal{T}\cong\mathbb{R}_{(2)}\mathcal{T}$.

\item If $W\in\Inv(P,\hbox{}^*\modfg_A)$, and $N$ is any graded $A$--module, then $\mathcal{T}(\roos^* (W),N)\cong\roos_* (\mathcal{T}(W,N))$, where $\roos_*$ (respectively, $\roos^*$) denotes the homological (respectively, cohomological) Roos complexes as defined for instance in \cite[\S 3]{AlvarezBoixZarzuelaspectralsequences}.

\item If $V\in\Inv(P,\hbox{}^*\modfg_A)$ is acyclic with respect to the limit, and $I$ is a graded, injective $A$--module, then there is a natural isomorphism of graded, $A$--modules
\[
\colim_{p\in P} \hbox{}^* \Hom_A (W_p,I)\cong\hbox{}^*\Hom_A \left(\lim_{p\in P}V_p,I\right).
\]
\end{enumerate}
    
\end{lm}

\begin{proof}
First, we prove part (i). By the own definition of $\mathcal{T}$, notice that $\mathcal{T}(-,E)$ is an exact functor for any injective, graded $A$--module $E$. In this way, we can apply \cite[Chapter V, Theorem 8.1]{CartanEilenberg1956} and we obtain that $\mathbb{R}\mathcal{T}$ and $\mathbb{R}_{(1)}\mathcal{T}$ are isomorphic. On the other hand, by the analogous for inverse systems of \cite[Lemma 6.3 and Proposition 7.1]{BrunBrunsRomer2007}, we know that any $W\in\Inv(P,\hbox{}^*\modfg_A)$ admits a projective resolution by inverse systems of the form
\[
\bigoplus_{i\in I_j} A(-t_{i,j})_{\leq q_{i,j}},
\]
where $I_j$ is a finite index set, $t_{i,j}\in\mathbb{Z}$ and $q_{i,j}\in P$. Since, for any $t\in\mathbb{Z}$ and any $q\in P$, $\mathcal{T}(A(-t)_{\leq q},-)$ is an exact functor, we can appeal again to \cite[Chapter V, Theorem 8.1]{CartanEilenberg1956} in order to guarantee that $\mathbb{R}\mathcal{T}$ and $\mathbb{R}_{(2)}\mathcal{T}$ are isomorphic. This concludes the proof of part (i).

Now, we prove part (ii). Let $W\in\Inv(P,\hbox{}^*\modfg_A)$, and let $N$ be any graded $A$--module. In order to check that $\mathcal{T}(\roos^* (W),N)\cong\roos_* (\mathcal{T}(W,N))$ what we plan to show is that, for any $k\geq 0$,
\[
\mathcal{T}(\roos^k (W),N)\cong\roos_k (\mathcal{T}(W,N))\text{ and }\mathcal{T}(d^k,N)\cong d_{k+1},
\]
where $\xymatrix@1{\roos^k (W)\ar[r]^-{d^k}& \roos^{k+1} (W)}$ is the $k$--coboundary map of the cohomological Roos complex, and, on the other hand, $\xymatrix@1{\roos_{k+1} (\mathcal{T}(W,N))\ar[rr]^-{d_{k+1}}& & \roos_k (\mathcal{T}(W,N))}$.

Indeed, fix $k\geq 0$. By the definition of the cohomological Roos complex,
\[
\roos^k (W)=\prod_{p_0<\ldots<p_k} W_{p_0\ldots p_k},
\]
where $W_{p_0\ldots p_k}=W_{p_0}$. Since $P$ is finite, we can canonically identify $\roos^k (W)$ with
\[
\bigoplus_{p_0<\ldots<p_k} W_{p_0\ldots p_k}.
\]
Now, since the Hom functor converts direct sums into direct products, we have that $\mathcal{T}(\roos^k (W),I)$ is
\[
\prod_{p_0<\ldots<p_k} \hbox{}^* \Hom_A (W_{p_0\ldots p_k},N).
\]
Again, by the finiteness of $P$, we can identify this direct product with the direct sum
\[
\bigoplus_{p_0<\ldots<p_k} \hbox{}^* \Hom_A (W_{p_0\ldots p_k},N),
\]
which is nothing but $\roos_k (\mathcal{T}(W,N))$. This shows that both complexes $\mathcal{T}(\roos^* (W),N)$ and $\roos_* (\mathcal{T}(W,N))$ have isomorphic pieces.

On the other hand, it is enough to check that the differentials agree on each direct summand (because the pieces are direct sums). Indeed, on $W_{p_0\ldots p_k}$, $d^k$ is defined as
\[
\phi_{p_1p_0}\circ\pi_{p_1\ldots p_k}+\sum_{l=1}^{k+1} (-1)^l \pi_{p_0\ldots\widehat{p_l}\ldots p_{k+1}},
\]
where $\xymatrix@1{W_{p_1}\ar[rr]^-{\phi_{p_1p_0}}& & W_{p_0}}$ is the structural map of the inverse system $W$ for $p_0<p_1$, and, on the other hand, $\xymatrix@1{\roos^k(W)\ar[rr]^-{\pi_{p_0\ldots p_k}}&  & W_{p_0\ldots p_k}}$ is the corresponding projection map. Now, because of the fact that the Hom functor is additive and $A$--linear, we have that $\hbox{}^* \Hom_A (d^k,N)$ is defined on $\hbox{}^* \Hom_A (W_{p_0\ldots p_k},N)$ by the formula
\[
\hbox{}^*\Hom_A (\pi_{p_1\ldots p_k},N)\circ\hbox{}^*\Hom_A (\phi_{p_1p_0},N)+\sum_{l=1}^{k+1}(-1)^l\cdot \hbox{}^*\Hom_A (\pi_{p_0\ldots\widehat{p_l}\ldots p_{k+1}},N).
\]
Moreover, we observe that, for any $p_0<\ldots<p_k$, $\hbox{}^*\Hom_A (\pi_{p_0\ldots p_k},N)$ is precisely the insertion map of the direct system $\mathcal{T}(W,N)$. This finally shows, by the definition of the homological Roos complex, that $\mathcal{T}(d^k,N)\cong d_{k+1}$. This concludes the proof of part (ii).

Finally, we move on proving part (iii). Let $V\in\Inv(P,\hbox{}^*\modfg_A)$ is acyclic with respect to the limit, and let $I$ be a graded, injective $A$--module. Since $V$ is acyclic with respect to the limit, we know that the coaugmented cohomological Roos complex $\xymatrix@1{0\ar[r]& \lim_{p\in P}V_p\ar[r]& \roos^* (V)}$ is exact. After applying $\mathcal{T}(-,I)$ we end up with the following acyclic complex, and notice that here we are using that $I$ is injective: $\xymatrix@1{\mathcal{T} (\roos^* (V),I)\ar[r]& \mathcal{T} (\lim_{p\in P}V,I)\ar[r]& 0.}$ Notice that
\[
\mathcal{T} \left(\lim_{p\in P}V,I\right)=\hbox{}^*\Hom_A \left(\lim_{p\in P}V_p,I\right).
\]
By part (ii), we know that $\mathcal{T}(\roos^* (V),I)\cong\roos_* (\mathcal{T}(V,I))$. But $\roos_* (\mathcal{T}(V,I))$ has as single non--zero homology group, namely $\colim_{p\in P} \hbox{}^* \Hom_A (V_p,I)$. This finally shows that there is a natural isomorphism
\[
\colim_{p\in P} \hbox{}^* \Hom_A (V_p,I)\cong\hbox{}^*\Hom_A \left(\lim_{p\in P}V_p,I\right).
\]
This concludes the proof of this result.
\end{proof}


Next result can be regarded as a corrected version of the statement of \cite[Theorem 6.5 and Corollary 6.6]{AlvarezBoixZarzuelaspectralsequences}. We reproduce the proof for the convenience of the reader because it also contains differences with respect to the one presented in \cite{AlvarezBoixZarzuelaspectralsequences}, as the reader will shortly appreciate. 

\begin{teo}\thlabel{sucesion espectral y colapso todo en uno: segunda parte}
Preserving the notations and assumptions established along \thref{otra construccion mas de lo mismo}, set $T:=\hbox{}^* \Hom_A (-,C)$ for some finitely generated, graded $A$--module $C$ with finite injective dimension, and let $V\in\Inv (P,\hbox{}^*\modfg_A)$ which is acyclic with respect to the limit functor. Then, the following assertions hold.

\begin{enumerate}[(i)]

\item There is a spectral sequence of the form
\[
E_2^{-i,j}=\mathbb{L}_i \colim_{p\in P} \mathbb{R}^j \mathcal{T} (V,C)\Longrightarrow \mathbb{R}^{j-i} T\left(\lim_{p\in P} V\right).
\]

\item In addition, assume that for any $p\in P$, $\R^j T \left(V_p\right)=0$ up to a single value of $j$ (namely, $h_p$), and that there is a canonical isomorphism of graded $A$--modules
\[
\mathbb{R}^j \mathcal{T} (V,C)\cong\bigoplus_{j=h_q} \left(\mathbb{R}^{h_q} T (V_q)\right)_q.
\]
Then, the previous spectral sequence can be expressed in the following way:
\[
E_2^{-i,j}=\bigoplus_{j=h_q} \R^{h_q} T \left(V_q\right)^{\oplus m_{i,q}}\xymatrix{ \ar@{=>}[r]_-{i}& }\R^{j-i}T \left(\lim_{p\in P}V_p\right),
\]
where $m_{i,q}:=\dim_{\K} \left(\widetilde{H}_{i-1} \left(\left(q,1_{\widehat{P}}\right);\K\right)\right)$. Moreover,
this spectral sequence degenerates at the $E_2$-sheet. Therefore, by \cite[5.2.5]{Weibel1994}, for each
$0\leq r\leq\dim (A)$ there is an increasing, finite filtration $\{H_k^r\}$ of $\R^r T(\lim_{p\in P}V_p)$ by graded
$A$-modules such that, for any $k\geq 0$,
\[
H_k^r/H_{k-1}^r\cong E_2^{-k,r+k}=\bigoplus_{\{q\in P\ \mid\ r+k=h_q\}} \R^{h_q} T \left(V_q\right)^{\oplus m_{k,q}},
\]
where we follow the convention that $H_{-1}^r =0$.

\end{enumerate}
\end{teo}

\begin{proof}
Choose $\xymatrix@1{0\ar[r]& C\ar[r]& I^0\ar[r]& \ldots}$ a finite injective resolution of $C$ by graded injective objects, which we know by assumption is finite. Now, after fixing $p\in P$ and applying the contravariant functor $T$ we obtain the following cochain complex of graded $A$--modules:
\[
\xymatrix{0\ar[r]& \hbox{}^* \Hom_A(V_p,C)\ar[r]& \hbox{}^*\Hom_A(V_p,I^0)\ar[r]& \hbox{}^*\Hom_A(V_p,I^1)\ar[r]& \ldots}
\]
It is easy to check that this collection of cochain complexes of graded $A$--modules yields the following cochain complex of direct systems of graded $A$--modules:
\[
\xymatrix{0\ar[r]& \mathcal{T} (V,C)\ar[r]& \mathcal{T}(V,I^0)\ar[r]& \mathcal{T} (V,I^1)\ar[r]& \ldots}
\]
Now, by applying the homological Roos complex \cite[\S 3.1]{AlvarezBoixZarzuelaspectralsequences} for computing the left derived functors of the colimit, we end up with the bicomplex $\roos^{-i}(\mathcal{T}(V,I^j))$, where the reader will easily note
that we put a minus in the index $i$ because we want to work with a bicomplex located in the second quadrant.

Moreover, we have to stress that the vertical differentials are the
ones of the homological Roos complex and the horizontal ones are the
induced by our above resolution of $C$; so, the bicomplex
$\roos_i (\mathcal{T}(V,I^j))$ produces two spectral
sequences; namely, the ones provided respectively by the first and
the second filtration of the previous bicomplex. In this way, the
first thing one should ensure is that both spectral sequences
converge and calculate their common abutment.

First, we want to check that our spectral sequences converge;
indeed, since $\roos^{-i} (\mathcal{T}(V,I^j))=0$ for all $i\gg 0$ (this is because we are working over a finite poset), we have that the bicomplex $\roos^{-i} (\mathcal{T}(V,I^j))$ has just a finite number of rows, which implies the convergence of both spectral sequences.

Now, notice that the $E_2$-page of one of these spectral
sequences is obtained by firstly computing the homology of the
columns and then computing the cohomology of the rows; indeed, the columns are exact, and therefore we only need to compute the cohomology of the row located in position $0$. Now, in this position, this row is precisely given by the complex whose spots are $\colim_{p\in P}\mathcal{T}(V,I^j)$ for each $j$. Now, by using part (iii) of \thref{key lemma} we can guarantee that the cohomology of this row can also be expressed as the cohomology of the cochain complex given by $\{(\hbox{}^*\Hom_A (\lim_{p\in P}V_p,I^j)\}_{j\geq 0}$. Keeping in mind this fact, we have, keeping in mind part (i) of \thref{key lemma} that the abutment is precisely $\R^* T \left(\lim_{p\in P} V\right)$.

Finally, the other spectral sequence that we can produce is the one obtained by firstly taking cohomology on the rows and then calculating the homology of the columns; in this case, one obtains as $E_1$-page $E_1^{-i,j}=\roos^{-i} \left(\R^j
\mathcal{T} (V,C)\right)$. In addition, since the boundary map of the $E_1$-page is the one of the Roos chain complex, and this chain complex computes the $i$th left derived functor of the colimit, its $E_2$-page turns out to be $E_2^{-i,j}=\Lf_i \colim_{p\in P} \R^j\mathcal{T} (V,C)$. Summing up, combining all the foregoing facts we obtain the spectral sequence
\[
E_2^{-i,j}=\mathbb{L}_i \colim_{p\in P} \mathbb{R}^j \mathcal{T} (V,C)\Longrightarrow \mathbb{R}^{j-i} T\left(\lim_{p\in P} V\right).
\]
in the category of graded $A$-modules; this proves part (i).

Now, let us prove part (ii). By our assumption, we have a canonical isomorphism of graded $A$--modules
\[
\mathbb{R}^j \mathcal{T} (V,C)\cong\bigoplus_{j=h_q} \left(\mathbb{R}^{h_q} T (V_q)\right)_q.
\]
Now, fix $i\in\mathbb{N}_0$. Applying to the above isomorphism the $i$th left derived functor of the colimit over $P$, we obtain the following canonical isomorphism of graded $A$--modules:
\[
\mathbb{L}_i \colim_{p\in P} \mathbb{R}^j \mathcal{T} (V,C)\cong\bigoplus_{j=h_q} \mathbb{L}_i \colim_{p\in P}\left(\mathbb{R}^{h_q} T (V_q)\right)_q.
\]
Now, by \cite[Lemma 4.10]{AlvarezBoixZarzuelaspectralsequences} we know that
\[
\mathbb{L}_i \colim_{p\in P}\left(\mathbb{R}^{h_q} T (V_q)\right)_q\cong\bigoplus\widetilde{H}_{i-1} ((q,1_{\widehat{P}}); \mathbb{R}^{h_q} T (V_q)).
\]
Now, by giving the same arguments explained along the proof of \cite[Theorem 4.11]{AlvarezBoixZarzuelaspectralsequences} we conclude that there is a canonical isomorphism of graded $A$--modules
\[
\mathbb{L}_i \colim_{p\in P}\left(\mathbb{R}^{h_q} T (V_q)\right)_q\cong\bigoplus_{j=h_q} \R^{h_q} T \left(V_q\right)^{\oplus m_{i,q}},
\]
just what we finally wanted to show.
\end{proof}




\section{On ultrametric functions}\label{section on ultrametric functions}

The goal of this section is to introduce the notion of the so--called ultrametric functions, mainly because the regularity of a finitely generated, graded module over a polynomial ring may be regarded as a particular example of ultrametric function. We want to mention here that this definition is not standard terminology and it is only coined for the purposes of this paper. The reader may like to compare our definition of ultrametric function with the notion of Sylvester matrix function studied in \cite{Virilirankfunctions}, and also with the approach developed in \cite[Lemma 2.1]{Santoni}.

\begin{df}\thlabel{introducing ultrametric functions}
Let $\mathcal{A}$ be an abelian category, let $G$ be a finitely generated, torsion free abelian group with $r=\operatorname{rank}(G)$, let $P(G)$ be the poset of all subsets of $G$, and let $\xymatrix@1{\mathcal{A}\ar[r]^-{\mu}& \widehat{M}}$ be a map, where
\[
\widehat{M}=\begin{cases} \mathbb{Z}\cup\{-\infty\},\text{ if }r=1,\\
P(G),\text{ if }r\geq 2.
\end{cases}
\]
We say that $\mu$ is an \textit{ultrametric function} if
\[
\mu (0)=\begin{cases} -\infty,\text{ if }r=1,\\
G,\text{ if }r\geq 2.
\end{cases}
\]
and, for any short exact sequence $\xymatrix@1{0\ar[r]& M'\ar[r]& M\ar[r]& M''\ar[r]& 0}$ in $\mathcal{A}$, we have
\[
\mu (M)\begin{cases} \leq\max\{\mu (M'),\ \mu (M'')\},\text{ if }r=1,\\
\supseteq \mu(M')\cap\mu (M''),\text{ if }r\geq 2.\end{cases}
\]
\end{df}
Here, we present some examples of ultrametric functions.

\begin{ex}\thlabel{examples of ultrametric functions}
The goal here is to introduce several examples of ultrametric functions that appear quite naturally in the realm of Commutative Algebra and Algebraic Geometry.

\begin{enumerate}[(i)]

\item Let
\[
A=\bigoplus_{n\geq 0} A_n
\]
be a positively graded ring such that $A=A_0 [A_1]$. For each integer $l\geq 0,$ one has, as consequence of \cite[16.2.15 (ii)]{BroSha}, that the Castelnuovo--Mumford regularity $\reg^l$ at and above level $l$ is an ultrametric function of $\hbox{}^*\modfg_A$.

\item Under the same assumptions of part (i), let $\mathcal{X}\subset\mathbb{N}$ be a set of indices. We have, by \cite[Lemma 2.3 (i)]{Cavigliatensorproducts} that the regularity $\reg^{\mathcal{X}}$ with respect to $\mathcal{X}$ defined by Caviglia in \cite[Definition 2.1]{Cavigliatensorproducts} is also an ultrametric function.

\item Let now $\mathcal{C}$ be the category of finitely generated modules over a Noetherian ring $R,$ and let $I\subset R$ be any ideal. Then, by \cite[Proposition 1.2.9]{BrunsHerzog1993} we know that $\mu=-\grade_I$ is also an ultrametric function.

\item Let now $\mathcal{B}$ be the category of finitely generated modules over a local ring $(R,\mathfrak{m})$ having finite projective dimension. Then, the projective dimension $\projdim$ is also an ultrametric function on $\mathcal{B}.$

\item More generally, let now $\mathcal{G}$ be the category of finitely generated modules over a commutative Noetherian ring having finite Gorenstein dimension \cite[(1.2.3)]{Gdimensionbook}. Then, the Gorenstein dimension $\Gdim$ is an ultrametric function on $\mathcal{G}$ because of \cite[(1.2.9)]{Gdimensionbook}. We recall that, for finitely generated modules with finite projective dimension, projective and Gorenstein dimension coincide, as it is shown in \cite[(1.2.10)]{Gdimensionbook}.

\item Let now $G$ be a finitely generated, torsion free abelian group, let $R$ be a polynomial ring over a field graded by $G,$ and let $\mathcal{C}$ be the category of finitely generated, graded $R$--modules. In this case, the Castelnuovo--Mumford regularity as defined by Maclagan and Smith in \cite{MaclaganSmith2004} is also an ultrametric function because of \cite[Lemma 7.1]{MaclaganSmith2004}.

\item Let $\mathcal{C}$ be the category of finitely generated $R$--modules, where $R$ is a commutative Noetherian ring, and let $I$ be an ideal of $R.$ The cohomological dimension $\cdim(I,-)$ with respect to $I$ is also an ultrametric function, as proved for instance by Divaani--Aazar, Naghipour and Tousi in \cite[Corollary 2.3 (i)]{DivaaniAazarNaghipourTousi2002}.

\item Let $R$ be a commutative Noetherian local ring of prime characteristic $p,$ let $F$ be the Frobenius map on $R,$ and denote by $\mathcal{C}$ the category of left $R[\Theta; F]$--modules. Denote by $\rho$ the $F$--depth (respectively, generalized $F$--depth) function of an object of $\mathcal{C}$, which have been studied for instance in the work by Maddox and Miller \cite{MaddoxMiller}. Then, thanks to \cite[Theorem 3.5]{MaddoxMiller}, one has that $\mu=-\rho$ is also an ultrametric function.

\end{enumerate}

\end{ex}
One of the main reasons for considering ultrametric functions is their behaviour with respect to filtrations; namely:

\begin{lm}\thlabel{ultrametric functions and filtrations}
Let $\mathcal{A}$ be an abelian category, let $\mu$ be an ultrametric function on $\mathcal{A},$ set $r:=\rank (G)$ as in \thref{introducing ultrametric functions}, and let $0=F_{-1}\subset F_0\subset F_1\subset\ldots\subset F_t=F$ be a filtration in $\mathcal{A}.$ Then, one has that, if $r=1,$ then
\[
\mu (F)\leq\max_{0\leq k\leq t}\left\{\mu\left(F_k/F_{k-1}\right)\right\}.
\]
On the other hand, if $r\geq 2$ then
\[
\mu (F)\supseteq\bigcap_{k=0}^t \mu\left(F_k/F_{k-1}\right).
\]
\end{lm}

\section{Main results}\label{section of main results}
As the title says, the main purpose of this section is to present the main results of this paper. Before that, we establish the following consequence of Theorem \ref{sucesion espectral y colapso todo en uno: segunda parte}.

\begin{teo}\thlabel{spectral sequence for the deficiency modules}
Let $\mathbb{K}$ be a field, and let $A=\mathbb{K}[x_1,\ldots ,x_n]$ be the polynomial ring in $n$ variables over $\mathbb{K}$ having the standard grading. Moreover, set $\omega_A:=A(-n),$ let $V$ be an inverse system of $\hbox{}^*\modfg_A$ indexed by a poset $P$ which is acyclic with respect to the limit on $P$. Finally, we denote by $\hbox{}^*\mathcal{H}om_A (-,-)$ the bivariate functor as defined in \thref{otra construccion mas de lo mismo}). Then, there exists a spectral sequence
\[
E_2^{-i,j}=\mathbb{L}_i \colim_{p\in P} \hbox{}^*\mathcal{E}xt_A^j (V_p, \omega_A)\xymatrix{\hbox{}\ar@{=>}[r]& }\hbox{}^*\Ext_A^{j-i} \left(\lim_{p\in P}V_p,\omega_A\right).
\]
\end{teo}

\begin{disc}\thlabel{the inequality for any ultrametric function}
Let $I$ be a homogeneous ideal of $A$, let $P$ be the poset given by the all the possible sums of the primary components of $I$ ordered by reverse inclusion, and let $A/[*]$ be the inverse system given, for each $p\in P$, by $A/I_p$, where the structural maps are either zero or the projections $A/I_p\longrightarrow A/I_q$ when $q<p$. Assume that the below statements hold.

\begin{enumerate}[(i)]

\item $P$ is a subset of a distributive lattice of ideals of $A$ such that $I+J\in P$ for all $I,J\in P$. Notice that this assumption implies \cite[Example 3.3]{BrunBrunsRomer2007} that the inverse system $A/[*]$ is acyclic with respect to the limit on $P$.

\item For any $p\in P,$ $A/I_p$ is a Cohen--Macaulay ring. Observe that this assumption implies that $\hbox{}^*\Ext_A^k (A/I_p,A(-n))$ is zero up to a single value of $k$.

\item If $p<q$ then $\height(I_q)<\height(I_p)$. This assumption ensures that there is a canonical isomorphism of graded $A$--modules
\[
\hbox{}^*\mathcal{E}xt_A^j(A/[*], A(-n))\cong\bigoplus_{j=h_q} \left(\hbox{}^*\Ext_A^{h_q} \left(A/I_q,A(-n)\right)\right)_q.
\]

\end{enumerate}
In this way, the spectral sequence obtained in \thref{spectral sequence for the deficiency modules} deduced from \thref{sucesion espectral y colapso todo en uno: segunda parte} boils down to the below one:
\[
E_2^{-i,j}=\bigoplus_{j=h_q} \hbox{}^*\Ext_A^{h_q} \left(V_q,A(-n)\right)^{\oplus m_{i,q}}\xymatrix{ \ar@{=>}[r]_-{i}& }\hbox{}^*\Ext_A^{j-i} \left(\lim_{p\in P}V_p,A(-n)\right),
\]
where $m_{i,q}:=\dim_{\K} \left(\widetilde{H}_{i-1} \left(\left(q,1_{\widehat{P}}\right);\K\right)\right)$.

Now, given a finitely generated, graded $A$--module, its $j$th deficiency graded module is given by
\[
K^j (M):=\hbox{}^*\Ext_A^{n-j} (M,A(-n)).
\]
Then, appealing again to \thref{sucesion espectral y colapso todo en uno: segunda parte}, we have that for any integer $j\geq 0$ there exists a filtration $\{H_k^{n-j}\}$ of
\[
K^j \left(\lim_{p\in P} A/I_p\right)
\]
such that, for each $k\geq 0,$
\[
H_k^{n-j}/H_{k-1}^{n-j}\cong\bigoplus_{\substack{p\in P\\ n-j+k=h_p}} K^{\dim (A/I_p)} (A/I_p)^{\oplus m_{n-j,p}},
\]
where we follow as usual the convention $H_{-1}^{n-j}=0.$ This is equivalent to say, since $n=h_p+\dim (A/I_p),$ that
\[
H_k^{n-j}/H_{k-1}^{n-j}\cong\bigoplus_{\substack{p\in P\\ j-k=\dim (A/I_p)}} K^{\dim (A/I_p)} (A/I_p)^{\oplus m_{n-j,p}},
\]
where
\[
m_{n-j,p}:=\dim_{\mathbb{K}}\left(\widetilde{H}_{j-\dim(A/I_p)-1}((p,1_{\widehat{P}});\mathbb{K})\right).
\]
In this way, given any ultrametric function $\mu$ on the category of finitely generated, $\mathbb{Z}$--graded $A$--modules we have, according to \thref{ultrametric functions and filtrations}, that
\begin{equation}\label{general inequality}
\mu\left(K^j \left(\lim_{p\in P} A/I_p\right)\right)\leq\max_{p\in S_j}\mu\left(K^{\dim (A/I_p)} (A/I_p)\right),
\end{equation}
where $S_j:=\{p\in P:\ j\geq\dim (A/I_p),\ m_{n-j,p}\neq 0\}$.
\end{disc}

Next statement is the main result of this paper, which recovers and extends Kummini--Murai's upper bound for the regularity of deficiency modules of monomial ideals.

\begin{teo}\thlabel{the expected bound when P is a C.M. poset}
Preserving the assumptions and notations fixed along \thref{the inequality for any ultrametric function}, we have that
\[
\reg\left(K^j\left(\lim_{p\in P} A/I_p\right)\right)\leq\max_{p\in S_j} (\dim (A/I_p))\leq j.
\]
\end{teo}

\begin{proof}
First of all, without loss of generality we can assume that $\depth(A/I)\leq j\leq\dim (A/I),$ otherwise $K^j (A/I)=0$ and the result is obviously true. As we observed in Example \ref{examples of ultrametric functions}, the regularity is an ultrametric function; therefore, according to \thref{the inequality for any ultrametric function} we have that
\[
\reg\left(K^j\left(\lim_{p\in P} A/I_p\right)\right)\leq\max_{p\in S_j}\reg(K^{\dim (A/I_p)}(A/I_p)).
\]
Moreover, since the poset $P$ satisfies, for any $p\in P,$ that $A/I_p$ is a Cohen--Macaulay ring and therefore, by applying \cite[Proposition 2.3]{HaHoa2008}, we obtain that
\[
\reg(K^{\dim (A/I_p)}(A/I_p))=\dim (A/I_p),
\]
just what we wanted to show.
\end{proof}

\begin{rk}
The reader will easily note that, given $j\geq 0$ an integer,
\[
\{p\in P:\ j=\dim (A/I_p),\ m_{n-j,p}\neq 0\}=\{p\in P:\ j=\dim (A/I_p)\}
\]
if and only if $(p,1_{\widehat{P}})=\emptyset;$ indeed, if $j=\dim (A/I_p)$ then
\[
\widetilde{H}_{j-\dim (A/I_p)-1}((p,1_{\widehat{P}});\mathbb{K})=\widetilde{H}_{-1}((p,1_{\widehat{P}});\mathbb{K}),
\]
and this homology group is not zero if and only if $(p,1_{\widehat{P}})=\emptyset,$ as claimed. From here and from \thref{the expected bound when P is a C.M. poset} we can deduce some interesting consequences.

\begin{enumerate}[(i)]

\item On the one hand, given $j\geq 0$, assume that there exists $p\in P$ such that $I_p$ is minimal and $\dim (A/I_p)=j$. This single fact implies that $K^j (A/I_p)\neq 0$.

\item On the other hand, in case $\reg (K^j (A/I))=j$, then we can guarantee that there exists $p\in P$ such that $\dim (A/I_p)=j$ and $m_{n-j,p}\neq 0$.

\end{enumerate}
\end{rk}

\begin{rk}
Let $\mathbb{K}$ be a field, let $A=\mathbb{K}[x_1,\ldots,x_n],$ graded in the standard way, and let $M$ be a finitely generated, graded $A$--module of dimension $d$. Borrowing terminology from \cite[Definition 2.3]{DaoMaVarbaro}, we say that $M$ satisfies the $r$--th Murai--Terai's condition $(MT_r)$ provided $\reg (K^i (M))\leq i-r$ for any $0\leq i\leq d-1$.

The reader will easily note that, under the assumptions of \thref{the inequality for any ultrametric function}, \thref{the expected bound when P is a C.M. poset} shows in particular that $A/I$ always satisfies $(MT_0)$. It is also shown in \cite[Propositions 3.1 and 4.1]{DaoMaVarbaro} that $A/I$ also satisfies $(MT_0)$ if either $A/I$ defines a Du Bois singularity in characteristic zero, or an $F$--pure singularity in prime characteristic.
\end{rk}

\section{The case of toric face rings}\label{section: the case of toric face rings}

The purpose of this section is to specialize \thref{the expected bound when P is a C.M. poset} in the case of toric face rings. Indeed, the following upper bound for the regularity of deficiency modules of certain toric face rings is, to the best of our knowledge, completely new. Once again, the following statement is an immediate consequence of \thref{the expected bound when P is a C.M. poset}. For unexplained terminology concerning toric face rings, we refer the reader, for instance, to Nguyen's paper \cite{HopNguyen2012}.

\begin{teo}\thlabel{deficiency in case of toric face rings}
Let $\Sigma\subseteq\mathbb{R}^d$ be a rational pointed fan, let $\mathcal{M}_{\Sigma}$ be a Cohen--Macaulay monoidal complex supported on $\Sigma,$ let $\mathbb{K}$ be a field, and let $B:=\mathbb{K}[\mathcal{M}_{\Sigma}]$ be the corresponding toric face ring. Moreover, we write $B:=A/I,$ where $A:=\mathbb{K}[x_1,\ldots, x_n]$ and $I$ is an ideal generated by squarefree monomials and pure binomials. Then, we have that
\[
\reg(K^j(B))\leq j.
\]
\end{teo}

\begin{proof}
We only need to check that conditions (i), (ii) and (iii) established along \thref{the inequality for any ultrametric function} hold. On the one hand, condition (i) was proved by Brun, Bruns and R\"{o}mer in \cite[Proposition 2.2 and Proof]{BrunBrunsRomer2007}. On the other hand, conditions (ii) and (iii) also hold because we are assuming that our monoidal complex is Cohen--Macaulay. The proof is therefore completed.
\end{proof}

\begin{rk}
The statement of \thref{deficiency in case of toric face rings} provides, to the best of our knowledge, a new bound for the regularity of these deficiency modules. Notice that our assumptions do not necessarily imply that our toric face ring is Cohen--Macaulay, see for instance \cite[Theorem 1.1]{IchimRomer2007} obtained by Ichim and R\"{o}mer.
\end{rk}

\section{The case of binomial edge ideals}\label{section: the case of binomial edge ideals}
The purpose of this section is to specialize our main results to obtain upper bounds of the regularity of the deficiency modules of the so--called binomial edge ideals, that were independently introduced in Herzog et al. \cite{binomialedgeideals} and Ohtani \cite{Ohtanibinomialedgeideals}, we refer to the interested reader to consult the recent book \cite[Chapter 7]{HerzogHibiOhsugibinomialidealsbook} for more information concerning this specific class of ideals. For our purposes here, we are interested in the following construction, which is borrowed from work by \`Alvarez Montaner (see \cite[Definition 3.3]{AlvarezMontaner2020}).

\begin{cons}\thlabel{Alvarez Montaner poset: construction}
Let $\mathbb{K}$ be a field, let $A=\mathbb{K}[x_1,\ldots,x_n,y_1,\ldots,y_n]$ be a polynomial ring over $\mathbb{K},$ and let $I\subseteq A$ be the binomial edge ideal associated to a graph $G$ in the set of vertices $[n].$ Moreover (see \cite[Theorem 7.16]{HerzogHibiOhsugibinomialidealsbook}), let
\[
J_G=P_{D_1} (G)\cap\ldots\cap P_{D_r} (G)
\]
be the minimal prime decomposition of $J_G$ (indeed, this is possible because $J_G$ is a radical ideal, see \cite[Corollary 7.13]{HerzogHibiOhsugibinomialidealsbook}) and let $P$ be the poset given by all the possible sums of the ideals in the decomposition ordered by reverse inclusion.

In this way, we construct in an iterative way a new poset $Q:=Q_{J_G}$ as follows: it is made up by the prime ideals contained in $P,$ and by the prime ideals appearing in the posets $P_{I_q}$ of sums of ideals in the decomposition of every non prime ideal $I_q$ of $P,$ and the prime ideals we obtain repeating this procedure every time we find a non prime ideal.
\end{cons}

Keeping in mind \thref{Alvarez Montaner poset: construction}, the main result of this section can be stated as follows.

\begin{teo}\thlabel{deficiency in case of binomial edge ideals}
Let $\mathbb{K}$ be a field, let $A=\mathbb{K}[x_1,\ldots,x_n,y_1,\ldots,y_n]$ be a polynomial ring over $\mathbb{K},$ and let $I\subseteq A$ be the binomial edge ideal associated to a graph $G$ in the set of vertices $[n]$. Assume that the inverse system $(A/I_q)_{q\in Q}$ is acyclic with respect to the limit functor. Then, we have that
\[
\reg\left(K^j\left(\lim_{q\in Q}A/I_q\right)\right)\leq j.
\]
\end{teo}

\begin{proof}
We only need to check that conditions (i), (ii) and (iii) established along \thref{the inequality for any ultrametric function} hold for the poset $Q$. Indeed, first of all condition (i) follows because we have supposed that $(A/I_q)_{q\in Q}$ is acyclic with respect to the limit functor. On the other hand, parts (ii) and (iii) are also true as observed along \cite[page 338, paragraph before Discussion 3.6]{AlvarezMontaner2020}. The proof is therefore completed.
\end{proof}

\begin{rk}\label{Alvarez Montaner remark}
There is a gap in the proof of \cite[Proposition 3.7]{AlvarezMontaner2020}, and thus we cannot guarantee the acyclicity with respect to the limit functor of the inverse system $(A/I_q)_{q\in Q}$. Is for this reason that we require this assumption in the statement of Theorem \ref{deficiency in case of binomial edge ideals}.
\end{rk}

\section{Examples}\label{section of examples}
The purpose of this final section is to illustrate several of our results by means of examples. In all these examples the unjustified claims and calculations were done with Macaulay2 \cite{M2}. Our first example illustrates that, in general, our bound is sharper for Kummini--Murai's one even in the monomial case.

\begin{ex}\thlabel{sometimes S_j is empty}
Let $I=(xz,xw,yz,yw)\subset R:=\mathbb{K}[x,y,z,w].$ In this case, we know that $I=(x,y)\cap (z,w),$ and that our poset has vertices
\[
(x,y)=p_1,\ (z,w)=p_2,\ p_3=(x,y,z,w).
\]
Moreover, the deficiency modules $K^j (R/I)$ are non--zero for $j=1$ and $j=2$. Moreover, we also have $(p_1,1_{\widehat{P}})=\emptyset=(p_2,1_{\widehat{P}}),\ (p_3,1_{\widehat{P}})=\{p_1\}\cup\{p_2\}$. This implies that
\begin{align*}
& \widetilde{H}_{j-\dim (R/I_{p_1})-1}((p_1,1_{\widehat{P}});\mathbb{K})\neq 0\Longleftrightarrow j=2,\\
& \widetilde{H}_{j-\dim (R/I_{p_2})-1}((p_2,1_{\widehat{P}});\mathbb{K})\neq 0\Longleftrightarrow j=2,\\
& \widetilde{H}_{j-\dim (R/I_{p_3})-1}((p_3,1_{\widehat{P}});\mathbb{K})\neq 0\Longleftrightarrow j=1.
\end{align*}
This shows that
\begin{align*}
& S_0=\{p\in P:\ 0\geq\dim (A/I_p),\ \widetilde{H}_{0-\dim(A/I_p)-1}((p,1_{\widehat{P}});\mathbb{K})\neq 0\}=\emptyset,\\
& S_1=\{p\in P:\ 1\geq\dim (A/I_p),\ \widetilde{H}_{1-\dim(A/I_p)-1}((p,1_{\widehat{P}});\mathbb{K})\neq 0\}=\{p_3\},\\
& S_2=\{p\in P:\ 2\geq\dim (A/I_p),\ \widetilde{H}_{2-\dim(A/I_p)-1}((p,1_{\widehat{P}});\mathbb{K})\neq 0\}=\{p_1, p_2\}.
\end{align*}
Summing up, we have that our bounds in this case are
\[
-\infty=\reg(K^0 (R/I))<0,\ 0=\reg (K^1 (R/I))\leq 0<1,\ 2=\reg (K^2 (R/I))\leq 2.
\]
This shows, in particular, that our bound for $\reg (K^1 (R/I))$ is better than Kummini--Murai's one, even in the monomial case.
\end{ex}
Our next example in the monomial case illustrates that in our upper bound for the regularity it is not only important the poset of sums, but also the dimensions of the ideals on it.

\begin{ex}\thlabel{non equdimensional example}
Let $I=(xy,xz)\subset R:=\mathbb{K}[x,y,z].$ In this case, we know that $I=(x)\cap (y,z),$ and that our poset has vertices
\[
(x)=p_1,\ (y,z)=p_2,\ (x,y,z)=p_3.
\]
Moreover, the deficiency modules $K^j (R/I)$ are non--zero for $j=1$ and $j=2.$ Indeed, $R/I$ is not Cohen--Macaulay because it is actually non--pure. Actually, in this case, whereas the simplicial complex attached to $I$ through the Stanley--Reisner correspondence is not Cohen--Macaulay, the order complex attached to our poset $P$ is Cohen--Macaulay. Now, it only misses to compute the sets $S_j$.

First of all, notice that
\begin{align*}
& \{p\in P:\ 0\geq\dim (R/I_p)\}=\{p_3\},\ \{p\in P:\ 1\geq\dim (A/I_p)\}=\{p_2,\ p_3\},\\
& \{p\in P:\ 2\geq\dim (A/I_p)\}=P.
\end{align*}
Moreover, we also have $(p_1,1_{\widehat{P}})=\emptyset=(p_2,1_{\widehat{P}}),\ (p_3,1_{\widehat{P}})=\{p_1\}\cup\{p_2\}$. This implies that
\begin{align*}
& \widetilde{H}_{j-\dim (R/I_{p_1})-1}((p_1,1_{\widehat{P}});\mathbb{K})\neq 0\Longleftrightarrow j=2,\\
& \widetilde{H}_{j-\dim (R/I_{p_2})-1}((p_2,1_{\widehat{P}});\mathbb{K})\neq 0\Longleftrightarrow j=2,\\
& \widetilde{H}_{j-\dim (R/I_{p_3})-1}((p_3,1_{\widehat{P}});\mathbb{K})\neq 0\Longleftrightarrow j=1.
\end{align*}
This shows that
\begin{align*}
& S_0=\{p\in P:\ 0\geq\dim (A/I_p),\ \widetilde{H}_{0-\dim(A/I_p)-1}((p,1_{\widehat{P}});\mathbb{K})\neq 0\}=\emptyset,\\
& S_1=\{p\in P:\ 1\geq\dim (A/I_p),\ \widetilde{H}_{1-\dim(A/I_p)-1}((p,1_{\widehat{P}});\mathbb{K})\neq 0\}=\{p_2, p_3\},\\
& S_2=\{p\in P:\ 2\geq\dim (A/I_p),\ \widetilde{H}_{2-\dim(A/I_p)-1}((p,1_{\widehat{P}});\mathbb{K})\neq 0\}=\{p_1\}.
\end{align*}
Summing up, we have that our bounds in this case are
\[
-\infty=\reg(K^0 (R/I))<0,\ 0=\reg (K^1 (R/I))=0,\ 2=\reg (K^2 (R/I))=2.
\]
\end{ex}

Now, we want to illustrate Theorem \ref{deficiency in case of toric face rings} by means of the following example, which is given by a non Cohen--Macaulay toric face ring borrowed from \cite[Example 4.2.3]{HopNguyenthesis}.

\begin{ex}\thlabel{a non C.M. toric face ring}
We consider the toric face ring
\[
R:=\frac{\mathbb{Q}[x_1,x_2,x_3,x_4,x_5,x_6]}{(x_1x_2-x_4^2,x_2x_3-x_5^2,x_4x_6,x_5x_6)}.
\]
In this case, this is a three dimensional ring and its unique non--zero deficiency modules are $K^2 (R)$ and $K^3 (R)$ with respective regularities $0$ and $3.$ Moreover, we have that
\[
I:=(x_1x_2-x_4^2,x_2x_3-x_5^2,x_4x_6,x_5x_6)=(x_2x_3-x_5^2,x_1x_2-x_4^2,x_6)\cap (x_1,x_3,x_4,x_5)\cap (x_2,x_4,x_5).
\]
is the prime decomposition of the defining ideal of $R.$ In this way, the poset $P$ of all the sums of the prime components of $I$ turns out to be
\begin{align*}
& p_1=(x_2x_3-x_5^2,x_1x_2-x_4^2,x_6),\ p_2=(x_1,x_3,x_4,x_5),\ p_3=(x_2,x_4,x_5),\\
& p_4=(x_1,x_3,x_4,x_5,x_6),\ p_5=(x_2,x_4,x_5,x_6),\ p_6=(x_1,x_2,x_3,x_4,x_5),\\
& p_7=(x_1,x_2,x_3,x_4,x_5,x_6),
\end{align*}
In this case, since
\begin{align*}
& S_2=\{p\in P:\ 2\geq\dim (A/I_p),\ \widetilde{H}_{2-\dim(A/I_p)-1}((p,1_{\widehat{P}});\mathbb{K})\neq 0\}=\{p_2,p_4,p_6,p_7\},\\
& S_3=\{p\in P:\ 3\geq\dim (A/I_p),\ \widetilde{H}_{3-\dim(A/I_p)-1}((p,1_{\widehat{P}});\mathbb{K})\neq 0\}=\{p_1,p_3,p_5\},
\end{align*}
we can guarantee, summing up, that our bound in this case is
\[
0=\reg (K^2 (R))\leq 2=2,\ 3=\reg (K^3(R))\leq 3.
\]
This example shows that in general the inequality
\[
\reg(K^j(A/I))\leq\max_{p\in S_j} (\dim (A/I_p))
\]
might be strict even in case $K^j (A/I)\neq 0$.
\end{ex}

We continue with an example of a binomial edge ideal borrowed from \cite[Example 3.2]{AlvarezMontaner2020}.

\begin{ex}
Let $A:=\mathbb{Q}[x_1,\ldots,x_5,y_1,\ldots,y_5],$ and let $J_G$ be the binomial edge ideal attached to the $5$--th path $G:=\{\{1,2\},\{2,3\},\{3,4\},\{4,5\}\}$. In this case, this is a six dimensional ring and its unique non--zero deficiency module $K^6 (R)$ has regularity $6.$ In this case, the poset $Q_{J_G}$ has $17$ elements and the regularities of the corresponding canonical modules of this list are (without repetitions) $\{3,4,5,6\}$.
\end{ex}

We finish with the below example, which is a non Cohen--Macaulay ring attached to a binomial edge ideal given by a complete bipartite graph. For more information about this specific kind of binomial edge ideals attached to complete bipartite graphs the reader can consult the work by Schenzel and Zafar, specially \cite[Section 4]{SchenzelZafarbinomialedgeideals}.

\begin{ex}
Let $A:=\mathbb{Q}[x_1,\ldots,x_8,y_1,\ldots,y_8],$ and let $J_G$ be the binomial edge ideal attached to the complete graph $K_{2,3}$. In this case, this is a ten dimensional ring and its unique non--zero deficiency modules are $K^j (R),$ where $j\in\{5,6,7,9,10\},$ with respective regularities $\{4,6,6,9,10\}$. In this case, the poset $Q_{J_G}$ has $6$ elements $\{p_1,\ldots, p_6\}$ such that
\[
\dim (A/I_{p_1})=10,\ \dim (A/I_{p_2})=9,\ \dim (A/I_{p_3})=6,\ \dim (A/I_{p_4})=6,\ \dim (A/I_{p_5})=0,\ \dim (A/I_{p_6})=4.
\]
The poset $Q_{J_G}$ in this specific example can be represented through the following Hasse--Voght diagram.
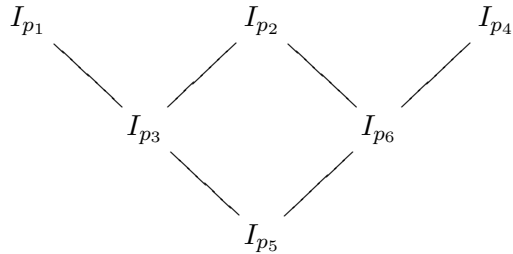
\[
\xymatrix{I_{p_1}&\hbox{}&  I_{p_2}& \hbox{}& I_{p_4}\\ \hbox{}& I_{p_3}\ar@{-}[ul]\ar@{-}[ur]& \hbox{}& I_{p_6}\ar@{-}[ul]\ar@{-}[ur]& \hbox{}\\ \hbox{}& \hbox{}&  I_{p_5}\ar@{-}[ul]\ar@{-}[ur]& \hbox{}& \hbox{}}
\]
Now, we want to write down the sets $S_j$'s for the non--zero canonical modules of $A/J_G$. We have
\begin{align*}
& S_5=\{p\in P:\ 5\geq\dim (A/I_p),\ \widetilde{H}_{5-\dim(A/I_p)-1}((p,1_{\widehat{P}});\mathbb{K})\neq 0\}=\{p_6\},\\
& S_6=\{p\in P:\ 6\geq\dim (A/I_p),\ \widetilde{H}_{6-\dim(A/I_p)-1}((p,1_{\widehat{P}});\mathbb{K})\neq 0\}=\{p_4\},\\
& S_7=\{p\in P:\ 7\geq\dim (A/I_p),\ \widetilde{H}_{7-\dim(A/I_p)-1}((p,1_{\widehat{P}});\mathbb{K})\neq 0\}=\{p_3\},\\
& S_9=\{p\in P:\ 9\geq\dim (A/I_p),\ \widetilde{H}_{9-\dim(A/I_p)-1}((p,1_{\widehat{P}});\mathbb{K})\neq 0\}=\{p_2\},\\
& S_{10}=\{p\in P:\ 10\geq\dim (A/I_p),\ \widetilde{H}_{10-\dim(A/I_p)-1}((p,1_{\widehat{P}});\mathbb{K})\neq 0\}=\{p_1\}.\\
\end{align*}
On the one hand, since $\dim (A/I_{p_6})=4$, we can guarantee that $4=\reg (K^5 (A/J_G))\leq 4<5$. On the other hand, since $\dim (A/I_{p_3})=6$, we can ensure that $6=\reg (K^7(A/J_G))\leq 6<7$. Hence once again our bound for this example turns out to be better than $\reg (K^j (A/J_G))\leq j$.
\end{ex}

\section*{Acknowledgements}
The authors would like to thank Josep \`Alvarez Montaner and Fernando Sancho de Salas for several interesting comments concerning the contents of this paper. We would also like to thank the referee for several useful comments and remarks that have substantially improved the contents and the correctness of this manuscript. Both authors received partial support by grant PID2022-137283NB-C22 funded by  MICIU/AEI/10.13039/501100011033.

\bibliographystyle{alpha}
\bibliography{AFBoixReferences}

\end{document}